\documentclass[12pt]{amsart}

\newtheorem{theorem}{Theorem}[section]

\newtheorem{lemma}[theorem]{Lemma}
\newtheorem{corollary}[theorem]{Corollary}
\theoremstyle{definition}   
\newtheorem{definition}{Definition}

\theoremstyle{remark}
\newtheorem{remark}[theorem]{Remark}

\numberwithin{equation}{section}
\usepackage[top=1in, bottom=1in, left=1in, right=1in]{geometry}
\usepackage{times}
\usepackage{enumerate}

\usepackage{graphicx,adjustbox}
\usepackage{hyperref}
\usepackage{amsmath,amssymb}
\usepackage{amscd}
\usepackage{graphicx}
\usepackage[all]{xy}

\title[Derivation modules and Poincar\'{e} series] 
{Derivation module and its Poincar\'{e} series of the 
projective closure of certain affine curves}
\author{
Joydip Saha
\and
Indranath Sengupta
\and
Pranjal Srivastava
}
\date{}
\address{\small \rm  Stat Math Unit, Indian Statistical Institute, Kolkata, West-Bengal,700108, INDIA.} 
\email{saha.joydip56@gmail.com}
\thanks{The first author thanks NBHM, Government of India for Post Doctoral Fellowship at ISI kolkata, through the research project PDF/2019/001074.}

\address{\small \rm  Discipline of Mathematics, IIT Gandhinagar, Palaj, Gandhinagar, 
Gujarat 382355, INDIA.}
\email{indranathsg@iitgn.ac.in}
\thanks{The second author is the corresponding author.}

\address{\small \rm  Discipline of Mathematics, IIT Gandhinagar, Palaj, Gandhinagar, 
Gujarat 382355, INDIA.}
\email{pranjal.srivastava@iitgn.ac.in}
\date{}

\subjclass[2020]{Primary 13N15, 13H10.}

\keywords{Monomial curves, module of derivation, Pseudo-Frobenius set, Poincar\'{e} series}

\begin{document}

\begin{abstract}
Our aim in this paper is to compute the Poincar\'{e} series of the derivation module of the 
projective closure of certain affine monomial curves. 
\end{abstract}

\maketitle

\section{Introduction}
Let $\mathbb{N}$ denote the set of nonnegative integers and $\mathbb{K}$ denote a field. 
Let $d\geq 3$, \, $\mathbf{\underline n} = (n_{1}, \ldots, n_{e})$ be a sequence of 
$e$ distinct positive integers with $\gcd(\mathbf{\underline n})=1$ and 
$n_{e}>n_{i}$ for all $i<e$. Let us assume 
that the numbers $n_{1}, \ldots, n_{e}$ generate the numerical semigroup 
$\Gamma(n_1,\ldots, n_e) = \left\lbrace\sum_{j=1}^{e}z_{j}n_{j}\mid z_{j}\in\mathbb{N}\right\rbrace$  
minimally, that is, if $n_i=\sum_{j=1}^{e}z_{j}n_{j}$ for some non-negative 
integers $z_{j}$, then $z_{j}=0$ for all $j\neq i$ and $z_{i}=1$. Let 
$\eta:\mathbb{K}[x_1,\,\ldots,\, x_e]\rightarrow \mathbb{K}[t]$ be the mapping defined by 
$\eta(x_i)=t^{n_i},\,1\leq i\leq e$. Let $\mathfrak{p} = \ker (\eta)$ be the   
defining ideal and 
$\mathbb{K}[\Gamma] = \mathbb{K}[x_1,\,\ldots,\, x_e]/\mathfrak{p}$ denote 
the semigroup ring. Let $n_{0}=0$ and define $\overline{\Gamma}\subset \mathbb{N}^{2}$ 
as the semigroup generated by the set of order pairs 
$\{(n_{i},n_{e}-n_{i})\mid 0\leq i\leq e\}$. 
Let $\overline{\mathfrak{p}(n_{1},\ldots n_{e})}$ denote the kernel of the 
$\mathbb{K}$-algebra map 
$\eta^{H}:\mathbb{K}[x_{0},\ldots,x_{e}]\longrightarrow \mathbb{K}[s,t]$, $\eta^{H}(x_{i})=t^{n_{i}}s^{n_{e}-n_{i}}, 0\leq i \leq e$. 
The homogenization of the ideal $\mathfrak{p}(n_1,\ldots, n_e) $ with respect to the variable $x_{0}$ is $\overline{\mathfrak{p}(n_{1},\ldots n_{e})}$. Let 
$\overline{C(n_{1},\ldots,n_{e})} = \{[(a^{n_{e}}:a^{n_{e}-n_{1}}b^{n_{1}}:\cdots:b^{n_{e}})]\in\mathbb{P}^{e}_{K}\mid a,b\in \mathbb{K}\}$, which is a projective curve, in fact 
the projective closure of the affine curve 
$C(n_{1},\ldots,n_{e}):=\{(b^{n_{1}},\ldots b^{n_{e}})\in \mathbb{A}^{e}_{\mathbb{K}}\mid b\in \mathbb{K}\}$. Let $\mathbb{K}[\overline{\Gamma}]$ denote the 
semigroup ring of the projective curve $\overline{C(n_{1},\ldots,n_{e})}$. 
The projective closure is called arithmetically Cohen-Macaulay if the 
semigroup ring $\mathbb{K}[\overline{\Gamma}]$ is Cohen-Macaulay, in other words 
the vanishing ideal 
$\overline{\mathfrak{p}(n_{1},\ldots n_{e})}$ is a Cohen-Macaulay ideal. 
Let $L=\{(a,b) \in \mathbb{Z}^{2} \mid a+b=n_{e}\mathbb{Z}\}$. For $i=1,2$, let 
$\Gamma_{i} = \Gamma_{i}(n_{1},\ldots,n_{e})$ denote the natural projection to the 
$i^{th}$ component of $\overline{\Gamma}$. Note that $\overline{\Gamma} \subset (\Gamma_{1} \times \Gamma_{2})\cap L$.
\medskip

Let $R$ be a commutative graded $\mathbb{K}$-algebra. 
The module of derivations 
$$\mathrm{Der}_{\mathbb{K}}(R)=\{\rho \in \mathrm{Hom}_{\mathbb{K}}(R,R) \vert \rho(ab)=a\rho(b)+b\rho(a) \,\forall \, a, b \in R\}$$ has an $R$-module structure. 
Let $\mu(\mathrm{Der}_{\mathbb{K}}(R))$ denote the minimal number 
of generators of $\mathrm{Der}_{\mathbb{K}}(R)$. In \cite{Kraft}, a 
set of generators of the module of 
derivation of numerical semigroup rings was computed using the pseudo-Frobenius set of 
the numerical semigroup. It has been observed 
that the minimal generator of $\mathrm{Der}_{k}(R)$ is related to the Cohen-Macaulay type 
by the relation $\mu(\mathrm{Der}_{\mathbb{K}}(\mathbb{K}[\Gamma]))= \rm{type}(\Gamma)+1$; 
see Patil et. al. \cite{Der-Arithmetic}. 
In the same paper, an explicit computation of a minimal generating set is given for 
the derivation module of an affine monomial curve associated to an almost arithmetic 
sequence. A generating set for the derivation module of the projective closure of 
monomial curves has been copmputed by Tamone et al. in \cite{Der}. 
\medskip

Finding a minimal free resolution of a finitely generated module $M$ over $R$ 
and their invariants, such as the Betti numbers are hard problems to solve. 
The Poincar\'{e} series of the finitely generated $R$-module $M$ is defined as 
$P_{M}^{R}(z)=\sum_{i\geq 0}\beta_{i}^{R}(M)z^{i}$, where $\beta_{i}^{R}(M)$ 
denotes the $i^{th}$ Betti number of $M$ over $R$. It is an interesting 
question to know whether the Poincar\'{e} series $P_{\mathbb{K}}^{R}(z)$ 
is a rational function; see 
\cite{Avramov}, \cite{Gulliksen}, \cite{Rossi}. This problem was 
studied by Fr\"{o}berg et al. \cite{Der affine}, where they computed  
the Poincar\'{e} series of the derivation module of a numerical semigroup ring 
$R = \mathbb{K}[\Gamma]$ and gave a relation between 
$P_{\mathrm{Der}_{\mathbb{K}}(R)}^{R}(z)$ 
and $P_{\mathbb{K}}^{R}(z)$. However, we are not 
aware of any such study on the projective closure $\mathbb{K}[\overline{\Gamma}]$ 
or on semigroup rings in higher dimension in general. 
In this article, we use Tamone's description from \cite{Der} to 
compute the Poincar\'{e} series of the derivation module of 
the coordinate ring $\mathbb{K}[\overline{\Gamma}]$ over $\mathbb{K}$, under the assumption 
that $\mathbb{K}[\overline{\Gamma}]$ is Cohen-Macaulay. We also show that the 
Poincar\'{e} series of the derivation module of the projective closure 
is related to the Poincar\'{e} series of the field $\mathbb{K}$ over 
$\mathbb{K}[\overline{\Gamma}]$ and the \emph{type} of $\Gamma_{i}$, for $i=1,2$. 
More generally, we prove that $P_{\mathbb{K}}^{\mathbb{K}[\overline{\Gamma}]}(z)$ 
is rational if and only if $P_{\mathrm{Der}_{\mathbb{K}}(\mathbb{K}[\overline{\Gamma}])}^{\mathbb{K}[\overline{\Gamma}]}(z)$ is rational. 
\medskip
 
We illustrate our results by studying the modules of derivations of the projective 
closure of two 
interesting classes of numerical semigroups defining monomial curves in 
$\mathbb{A}_{\mathbb{K}}^{4}$, viz., 
the Arslan curve and the Backelin curve. Feza Arslan introduced the class 
of curves defined by the numerical semigroup  
$\Gamma_{A_{h}}=\langle h(h+1),h(h+1)+1,(h+1)^{2},(h+1)^{2}+1 \rangle$, for $h\geq 2$, 
and proved the Cohen-Macaulayness of the tangent cone of the curve. Later, Herzog and 
Stamate \cite{CMC} studied the Cohen-Macaulayness of the projective closure of 
the Arslan curve. Subsequent to this, Fr\"{o}berg et al. presented the numerical 
semigroup $\Gamma_{B_{nr}}=\langle s,s+3,s+3n+1,s+3n+2 \rangle$, 
for $ n \geq 2, r\geq 3n+2, s=r(3n+2)+3$ 
and proved the unboundedness of the $\rm{type}(\Gamma_{B_{nr}})$. They 
attributed the claim to Backelin (See \cite{Froberg-Backelin}). Both the families 
of curves have the similarity that their projective closure are arithmetically 
Cohen-Macaulay (see \cite{CMC, Backelin}) and their Betti 
numbers are unbounded.

\section{Poincar\'{e} series}
In this section we prove Theorem  \ref{Poincare-1}, which allows us to determine 
the Poincar\'{e} series of $\mathrm{Der}_{\mathbb{K}}(R)$ over $R$, in terms of 
the Poincar\'{e} series of $\mathbb{K}$ over $R$ and the types of the numerical semigroups 
$\mathrm{pr}_{1}(\overline{\Gamma})=\Gamma_{1}=\langle n_{1},\dots,n_{e} \rangle$ 
and $\mathrm{pr}_{2}(\overline{\Gamma})=\Gamma_{2}=\langle n_{e}-n_{1},\dots,n_{e}-n_{e-1},n_{e} \rangle$. 
Let us first recall the theorem proved in \cite{Der}, which gives the structure 
of the derivation module. 

\begin{theorem}[\cite{Der}] \label{generator}
Let $R=\mathbb{K}[u^{n_{e}},v^{n_{1}}u^{n_{e}-n_{1}},\dots,v^{n_{e-1}}u^{n_{e}-n_{e-1}},v^{n_{e}}] $ be a Cohen-Macaulay ring. Then $\mathrm{Der}_{\mathbb{K}}(R)$ is generated by 
$D_{1} \cup \{u \frac{\partial} {\partial u}\} \cup D_{2} \cup \{v \frac{\partial} {\partial v}\}$, where 
$D_{1}$ and $D_{2}$ are as described below.
\begin{enumerate}

\item[1.] 
\begin{itemize}
\item[(i)] If $\Gamma_{2} \neq \mathbb{N}$, then $D_{1}=\{D_{\alpha}\vert \alpha-1 \in \mathrm{PF}(\Gamma_{2}\}=\{{v^{\beta}u^{\alpha}\frac{\partial} {\partial u}} \vert \alpha-1 \in \mathrm{PF}(\Gamma_{2})\}$, where 
$\beta$ is the least positive integer such that $(\beta,\alpha-1)\in (\Gamma_{1}\times \mathrm{PF}(\Gamma_{2}))\cap L$ and $(\beta,\alpha-1) +(n,n_{e}-n)\in \overline{\Gamma}$, for each $n \in \{0,n_{1},\dots,n_{e-1}\}$. 

\item[(ii)] If $\Gamma_{2} = \mathbb{N}$, then 
$D_{1} = \{D_{0}=v^{1+c'n_{e}}\frac{\partial} {\partial u}\}$, where $c'$ is the 
least natural number such that $(1+c'n_{e},-1) +(n,n_{e}-n)\in \overline{\Gamma}$, 
for each $n \in \{0,n_{1},\dots,n_{e-1}\}$.
\end{itemize}
\medskip

\item[2.]
\begin{itemize}
\item[(i)] If $\Gamma_{1} \neq \mathbb{N}$, then 
$D_{2} = \{D_{\delta}\vert \delta-1 \in \mathrm{PF}(\Gamma_{1})\}=\{{v^{\delta}u^{\gamma}\frac{\partial} {\partial v}} \vert \delta-1 \in \mathrm{PF}(\Gamma_{1})\}$, such that  $\gamma$ is the least positive integer 
with the properties $(\delta-1 ,\gamma)\in (\mathrm{PF}(\Gamma_{1})\times \Gamma_{2})\cap L$ and 
for each $n \in \{n_{1},\dots,n_{e}\}$, $\, \mathbf{(\delta-1 ,\gamma) +(n,n_{e}-n)\in \overline{\Gamma}}. \hfill (*)$

\item[(ii)] If $\Gamma_{1} = \mathbb{N}$, then 
$D_{2} =  \{D_{0}'=u^{1+e'n_{e}}\frac{\partial} {\partial v}\}$, where $e'$ is the least natural 
number such that $(-1,1+e'n_{e}) +(n,n_{e}-n)\in \overline{\Gamma}$, for each 
$n \in \{n_{1},\dots,n_{e-1}\}$.
\end{itemize}
\end{enumerate}
\end{theorem}

\begin{lemma}\label{ideal}
Let 
$R=\mathbb{K}[u^{n_{e}},v^{n_{1}}u^{n_{e}-n_{1}},\dots,v^{n_{e-1}}u^{n_{e}-n_{e-1}},v^{n_{e}}]$ 
be a Cohen-Macaulay ring. Suppose $\Gamma_{1}=\langle n_{1},\dots,n_{e} \rangle$ 
is the numerical semigroup with 
$\mathrm{PF}(\Gamma_{1})=\{\delta_{1}-1,\dots,\delta_{h_{1}}-1\}$ 
and  $\Gamma_{2}=\langle n_{e}-n_{1},\dots,n_{e}-n_{e-1},n_{e} \rangle$ 
is the numerical semigroup with 
$\mathrm{PF}(\Gamma_{2})=\{\alpha_{1}-1,\dots,\alpha_{h_{2}}-1\}$. 
Then the following statements are true:
\medskip

\begin{enumerate}
\item[(i)]
\begin{itemize}
\item[(a)] If $\Gamma_{2}\neq \mathbb{N}$, then $\bar{D}_{1}=(u,v^{\beta_{1}}u^{\alpha_{1}},\dots,v^{\beta_{h}}u^{\alpha_{h_{2}}})$ 
is isomorphic as a $\mathbb{K}[\overline{\Gamma}]$-module to the ideal 
$I=(u^{n_{e}},v^{\beta_{1}}u^{\alpha_{1}-1+n_{e}},\dots,v^{\beta_{h_{2}}}u^{\alpha_{h_{2}}-1+n_{e}})$.
\medskip

\item[(b)] If $\Gamma_{2}= \mathbb{N}$, then $\bar{D}_{1}=(u,v^{1+c'n_{e}})$ is isomorphic as a $\mathbb{K}[\overline{\Gamma}]$-module to the ideal $I=(u^{n_{e}},v^{1+c'n_{e}}u^{n_{e}-1})$. 
\end{itemize}
\medskip

\item[(ii)]
\begin{itemize}
\item[(a)] If  $\Gamma_{1}\neq \mathbb{N}$, then 
$\bar{D}_{2}=(v,u^{\gamma_{1}}v^{\delta_{1}},\dots,v^{\delta_{h_{1}}}u^{\gamma_{h_{1}}})$ is isomorphic as a $\mathbb{K}[\overline{\Gamma}]$-module to the ideal $J=(v^{n_{e}},v^{\delta_{1}+n_{e}-1}u^{\gamma_{1}},\dots,v^{\delta_{h_{1}}+n_{e}-1}u^{\gamma_{h_{1}}})$.
\medskip

\item[(b)] If $\Gamma_{1}= \mathbb{N}$, then $\bar{D}_{2}=(v,u^{1+e'n_{e}})$ 
is isomorphic as a $\mathbb{K}[\overline{\Gamma}]$-module to the ideal  $I=(v^{n_{e}},u^{1+e'n_{e}}v^{n_{e}-1})$.
\end{itemize}
\end{enumerate} 
\end{lemma}
\noindent (Note that $(\beta_{i},\alpha_{i}-1),(\delta_{i}-1,\gamma_{i}),(-1,1+e'n_{s})$ 
and $(1+c'n_{e},-1)$ 
satisfy the conditions of Lemma \ref{generator}.)

\proof
\begin{enumerate}[(i)]

\item (a) For $\Gamma_{2}\neq \mathbb{N}$, we can write 
\begin{align*}
(u,v^{\beta_{1}}u^{\alpha_{1}},\dots,v^{\beta_{h}}u^{\alpha_{h_{2}}})\cong &\, u^{n_{e}}(u,v^{\beta_{1}}u^{\alpha_{1}},\dots,v^{\beta_{h_{2}}}u^{\alpha_{h_{2}}})\\
=&\,(u^{n_{e}+1},v^{\beta_{1}}u^{\alpha_{1}+n_{e}},\dots,v^{\beta_{h_{2}}}u^{\alpha_{h_{2}}+n_{e}})\\
=& \,u(u^{n_{e}},v^{\beta_{1}}u^{\alpha_{1}-1+n_{e}},\dots,v^{\beta_{h_{2}}}u^{\alpha_{h_{2}}-1+n_{e}})\\
\cong & \,(u^{n_{e}},v^{\beta_{1}}u^{\alpha_{1}-1+n_{e}},\dots,v^{\beta_{h_{2}}}u^{\alpha_{h_{2}}-1+n_{e}}).
\end{align*}
Moreover, $(\beta_{i},\alpha_{i}-1)+(0,n_{e}) \in \overline{\Gamma}$, 
for each $i \in \{1,\dots,h_{2}\}$, therefore 
$$\bar{D}_{1}=(u,v^{\beta_{1}}u^{\alpha_{1}},\dots,v^{\beta_{h}}u^{\alpha_{h_{2}}})\cong (u^{n_{e}},v^{\beta_{1}}u^{\alpha_{1}-1+n_{e}},\dots,v^{\beta_{h_{2}}}u^{\alpha_{h_{2}}-1+n_{e}})$$ as $\mathbb{K}[\overline{\Gamma}]$-modules.
\medskip

(b) For $\Gamma_{2} = \mathbb{N}$, we can write
\begin{align*}
(u,v^{1+c'n_{e}})\cong &u^{n_{e}}(u,v^{1+c'n_{e}})\\
=&(u^{n_{e}+1},u^{n_{e}}v^{1+c'n_{e}})\\
=& u(u^{n_{e}},u^{n_{e}-1}v^{1+c'n_{e}})
\cong  (u^{n_{e}},u^{n_{e}-1}v^{1+c'n_{e}}).
\end{align*}
Moreover, $(1+c'n_{e},-1) +(0,n_{e})\in \overline{\Gamma}$, this implies 
$(u^{n_{e}},u^{n_{e}-1}v^{1+c'n_{e}}) \subset \mathbb{K}[\overline{\Gamma}]$ and 
$(u,v^{1+c'n_{e}})$ is isomorphic to the ideal $(u^{n_{e}},u^{n_{e}-1}v^{1+c'n_{e}})$ 
as $\mathbb{K}[\overline{\Gamma}]$-modules.
\medskip

\item (a)
When $\Gamma_{1}\neq \mathbb{N}$,
we have
\begin{align*}
(v,u^{\gamma_{1}}v^{\delta_{1}},\dots,v^{\delta_{h_{1}}}u^{\gamma_{h_{1}}})\cong & \,v^{n_{e}}(v,v^{\delta_{1}}u^{\gamma_{1}},\dots,v^{\delta_{h_{1}}}u^{\gamma_{h_{1}}})\\
=&(v^{n_{e}+1},v^{\delta_{1}+n_{e}}u^{\gamma_{1}},\dots,v^{\delta_{h_{1}}+n_{e}}u^{\gamma_{h_{1}}})\\
=& v(v^{n_{e}},v^{\delta_{1}-1+n_{e}}u^{\gamma_{1}},\dots,v^{\delta_{h_{1}}-1+n_{e}}u^{\gamma_{h_{1}}})\\
\cong & \,(v^{n_{e}},v^{\delta_{1}-1+n_{e}}u^{\gamma_{1}},\dots,v^{\delta_{h_{1}}-1+n_{e}}u^{\gamma_{h_{1}}}). 
\end{align*}
Since $(\delta_{i}-1,\gamma_{i})+(n_{e},0) \in \overline{\Gamma}$, for each $i \in \{1,\dots,h_{1}\}$, therefore  $$\bar{D}_{2}=(v,u^{\gamma_{1}}v^{\delta_{1}},\dots,v^{\delta_{h_{1}}}u^{\gamma_{h_{1}}}) \cong (v^{n_{e}},v^{\delta_{1}+n_{e}-1}u^{\gamma_{1}},\dots,v^{\delta_{h_{1}}+n_{e}-1}u^{\gamma_{h_{1}}})$$ as $\mathbb{K}[\overline{\Gamma}]$-modules.
\medskip

(b) Similarly for $\Gamma_{1}= \mathbb{N}$, we can write
\begin{align*}
(v,u^{1+e'n_{e}})\cong &v^{n_{e}}(v,u^{1+e'n_{e}})\\
=&(v^{n_{e}+1},u^{1+e'n_{e}}v^{n_{e}})\\
=&v(v^{n_{e}},u^{1+e'n_{e}}v^{n_{e}-1})\cong (v^{n_{e}},v^{n_{e}-1}u^{1+e'n_{e}}).
\end{align*}
Moreover, $(-1,1+e'n_{e}) +(n_{e},0)\in \overline{\Gamma}$, this implies $(v^{n_{e}},v^{n_{e}-1}u^{1+e'n_{e}}) \subset \mathbb{K}[\overline{\Gamma}]$ and $(v,u^{1+e'n_{e}})$ is isomorphic 
to the ideal $(v^{n_{e}},v^{n_{e}-1}u^{1+e'n_{e}})$ as $\mathbb{K}[\overline{\Gamma}]$-modules.
\end{enumerate}

\begin{lemma}\label{iso k}
\begin{enumerate}[(i)]
 \item $(\bar{v}^{\delta_{i}-1+n_{e}}\bar{u}^{\gamma_{i}})\bar{v}^{n_{i}}\bar{u}^{n_{e}-n_{i}}=\bar{0}$ in $\bar{R}=\frac{R}{(u^{n_{e}},v^{n_{e}})}$.

\item $(\bar{v}^{\beta_{i}}\bar{u}^{\alpha_{i}-1+n_{e}})\bar{v}^{n_{i}}\bar{u}^{n_{e}-n_{i}}=\bar{0}$ in $\bar{R}=\frac{R}{(u^{n_{e}},v^{n_{e}})}$. 
\end{enumerate}
Moreover,  $\bar{J}=(\bar{v}^{\delta_{1}+n_{e}-1}\bar{u}^{\gamma_{1}},\dots,\bar{v}^{\delta_{h_{1}}+n_{e}-1}\bar{u}^{\gamma_{h_{1}}})\cong \mathbb{K} \oplus \mathbb{K} 
\dots\oplus \mathbb{K} \, (h_{1}\textrm{-times})$ and \\$\bar{I}=(\bar{v}^{\beta_{1}}\bar{u}^{\alpha_{1}-1+n_{e}},\dots,\bar{v}^{\beta_{h_{2}}}\bar{u}^{\alpha_{h_{2}}-1+n_{e}})\cong \mathbb{K}\oplus \mathbb{K} \oplus \dots \oplus \mathbb{K}$ ($h_{2}$-times).
\end{lemma}

\proof Consider, 
\begin{align*}
(\bar{v}^{\delta_{i}-1+n_{e}}\bar{u}^{\gamma_{i}})\bar{v}^{n_{i}}\bar{u}^{n_{e}-n_{i}}=&\bar{v}^{\delta_{i}-1+n_{e}+n_{i}}\bar{u}^{(\gamma_{i}+n_{e}-n_{i})}\\
=&\bar{v}^{\delta_{i}-1+n_{i}}v^{n_{e}}\bar{u}^{(\gamma_{i}+n_{e}-n_{i})}((\delta_{i}-1,\gamma_{i}) +(n_{i},n_{e}-n_{i})+(n_{e},0)\in \overline{\Gamma})\\
=&\bar{0}.
\end{align*} 
\allowdisplaybreaks
Thus, the ideal $(\bar{v}^{\delta_{i}-1+n_{e}}\bar{u}^{\gamma_{i}})$ 
is isomorphic to $\mathbb{K}$ for each $i \in \{1,\dots,h_{1}\}$, 
hence  $\bar{J}\cong \mathbb{K}\oplus \mathbb{K} \oplus \dots \oplus \mathbb{K}$ 
($h_{1}$-times). Similarly $\bar{I}\cong \mathbb{K} \oplus \mathbb{K} \oplus \dots\oplus \mathbb{K} \, (h_{2}\textrm{-times})$. By the same argument as above we find that 
\begin{align*}
(\bar{v}^{n_{e}-1}\bar{u}^{1+e'n_{e}})\bar{u}^{n_{e}-n_{i}}\bar{v}^{n_{i}}= & 
\bar{v}^{n_{i}+n_{e}-1}\bar{u}^{1+e'n_{e}+n_{e}-n_{i}}\\
= & \bar{v}^{n_{i}+n_{e}-1}\bar{u}^{1+e'n_{e}+n_{e}-n_{i}}\\
= & \bar{v}^{n_{e}}\bar{v}^{-1+n_{i}}\bar{u}^{n_{e}-n_{i}+1+e'n_{e}}=\bar{0},
\end{align*}
since $(-1,1+e'n_{e}) +(n,n_{e}-n)\in \overline{\Gamma}$ for each 
$n \in \{n_{1},\dots,n_{e-1}\}$. Hence, the ideal 
$(\bar{v}^{n_{e}-1}\bar{u}^{1+e'n_{e}}) \cong \mathbb{K}$. 
Similarly, $(\bar{u}^{n_{e}},\bar{u}^{n_{e}-1}\bar{v}^{1+c'n_{e}})\bar{u}^{n_{e}-n_{i}}\bar{v}^{n_{i}}=\bar{0}$ implies that the ideal 
$(\bar{v}^{n_{e}-1}\bar{u}^{1+e'n_{e}}) \cong \mathbb{K}$. \qed

\begin{definition}
Suppose $U$ is an $R$-module. The Poincar\'{e} series of $U$ over R is 
$P_{U}^{R}(z)=\sum_{i \geq 0}\beta_{i}^{R}(U)z^{i}$, 
where $\beta_{i}$ is the $i^{th}$- total Betti number of $U$ over $R$.
\end{definition}

\begin{definition}(\cite{large})\rm{
If $A$ and $B$ are graded rings and $f:A \rightarrow B$ is a graded homomorphism which is surjective then $f$ is said to be large if $f_{*}:\mathrm{Tor}^{A}(\mathbb{K},\mathbb{K}) \rightarrow \mathrm{Tor}^{B}(\mathbb{K},\mathbb{K})$ is surjective.}
\end{definition}

\begin{lemma}[Theorem 1.1 \cite{large}]\label{Poincare}
Let $ (R,m,\mathbb{K})$ and $(S,n,\mathbb{K})$ are graded rings and $f:R \rightarrow S$ a graded homomorphism which is surjective. Then the following are equivalent
\begin{enumerate}
\item[(i)] The homomorphism $f$ is large, i.e. $f_{\star}: \mathrm{Tor}^{R}(\mathbb{K},\mathbb{K}) \rightarrow \mathrm{Tor}^{S}(\mathbb{K},\mathbb{K})$ is surjective.

\item[(ii)] For any finitely generated $S-$module $M$, $$P_{R}^{M}(z)=P_{S}^{M}(z)P_{R}^{S}(z).$$
\end{enumerate}
\end{lemma}

\begin{lemma}\label{COM P}
Let 
$R=\mathbb{K}[u^{n_{e}},v^{n_{1}}u^{n_{e}-n_{1}},\dots,v^{n_{e-1}}u^{n_{e}-n_{e-1}},v^{n_{e}}]$ 
be a Cohen-Macaulay ring and $M$ a finitely generated $R-$module such that 
$(u^{n_{e}},v^{n_{e}})M=0.$ Then $P_{M}^{R}(z)=P_{M}^{\bar{R}}(z)P_{\bar{R}}^{R}(z)=(1+z)P_{M}^{\bar{R}}(z).$
\end{lemma}

\proof Since $R$ is Cohen-Macaulay, $u^{n_{e}}$ is a non zero-divisor. By Theorem 2.2 in \cite{large}, $f_{1}:R \rightarrow \frac{R}{u^{n_{e}}}$ is large. Similarly, $v^{n_{e}}$ 
is a non zero-divisor in $ \frac{R}{u^{n_{e}}}$ and $f_{2}: \frac{R}{u^{n_{e}}} \rightarrow \frac{R}{(u^{n_{e}},v^{n_{e}})}$ is large. The composition of surjective maps is surjective, 
hence $f:R \rightarrow \frac{R}{(u^{n_{e}},v^{n_{e}})}$ is large. 
By Lemma \ref{Poincare}, $P_{M}^{R}(z)=P_{M}^{\bar{R}}(z)P_{\bar{R}}^{R}(z)=(1+z)P_{M}^{\bar{R}}(z).$
\medskip

\begin{theorem}\label{Poincare-1}
Let $R=\mathbb{K}[v^{n_{e}},u^{n_{1}}v^{n_{e}-n_{1}},\dots,u^{n_{e-1}}v^{n_{e}-n_{e-1}},u^{n_{e}}]$ be a Cohen-Macaulay ring. Let $\Gamma_{1}=\langle n_{1},\dots,n_{e} \rangle $ and 
$\Gamma_{2}=\langle n_{e}-n_{1},\dots,n_{e}-n_{e-1},n_{e} \rangle$ 
be the numerical semigroups with $\mathrm{PF}(\Gamma_{1})=\{\delta_{1}-1,\dots,\delta_{h_{1}}-1\}$ 
(if $\Gamma_{1} \neq \mathbb{N}$) and $\mathrm{PF}(\Gamma_{2})=\{\alpha_{1}-1,\dots,\alpha_{h_{2}}-1\}$ (if $\Gamma_{2} \neq \mathbb{N}$). Then 
\begin{itemize}
\item[(i)] 
\[
    P_{\bar{D}_{1}}^{R}(z)= 
\begin{cases}
    1+h_{2}P_{\mathbb{K}}^{R}(z),& \text{if } \Gamma_{2} \neq \mathbb{N}\\
    1+P_{\mathbb{K}}^{R}(z),              & \text{otherwise}
\end{cases}
\]

\item[(ii)] \[
    P_{\bar{D}_{1}}^{R}(z)= 
\begin{cases}
   1+h_{1}P_{\mathbb{K}}^{R}(z),& \text{if } \Gamma_{1} \neq \mathbb{N}\\
    1+P_{\mathbb{K}}^{R}(z),              & \text{otherwise}
\end{cases}
\]

\item[(iii)] 
    $P_{\mathrm{Der}_{\mathbb{K}}(R)}(z)= 1+(h_{1}+h_{2})P_{\mathbb{K}}^{R}(z)$.

\end{itemize}
Moreover, $P_{\mathbb{K}}^{R}(z)$ is rational if and only if $P_{\mathrm{Der}_{\mathbb{K}}(R)}^{R}(z)$ is rational.

(Note that $(\beta_{i},\alpha_{i}-1),(\delta_{i}-1,\gamma_{i}),(-1,1+e'n_{e})$ and $(1+c'n_{e},-1)$ satisfies the conditions of Lemma \ref{generator} respectively.)
\end{theorem}

\proof 
\begin{enumerate}[(i)]
\item By Lemma \ref{ideal}, when  $\Gamma_{2} \neq \mathbb{N}$, $\bar{D}_{1} \cong I= (u^{n_{e}},v^{\beta_{1}}u^{\alpha_{1}-1+n_{e}},\dots,v^{\beta_{h_{2}}}u^{\alpha_{h_{2}}-1+n_{e}})$ in $R$. Note that $P^{R}_{\bar{D}_{1}}(z)=P_{I}^{R}(z).$ 
 Suppose that $\bar{I}=(\bar{v}^{\beta_{1}}\bar{u}^{\alpha_{1}-1+n_{e}},\dots,\bar{v}^{\beta_{h_{2}}}\bar{u}^{\alpha_{h_{2}}-1+n_{e}})$ is the image of $I$ in $\bar{R}=\frac{R}{(u^{n_{e}},v^{n_{e}})}$, which is isomorphic to 
 $\mathbb{K}\oplus \mathbb{K} \oplus \dots \oplus \mathbb{K}$ 
 ($h_{2}$ copies), by Lemma \ref{iso k}. Since $R/I$ and $\bar{R}/\bar{I}$ are 
isomorphic as $\bar{R}$-modules, this implies 
$P^{\bar{R}}_{R/I}(z)=P^{\bar{R}}_{\bar{R}/\bar{I}}(z)$. 
Using Lemma \ref{COM P}, with $M=R/I$, and Lemma \ref{iso k}, we have
\begin{align*}
P_{I}^{R}(z)=&\frac{P_{R/I}(z)-1}{z}\\
=& \frac{(1+z)P^{\bar{R}}_{R/I}(z)-1}{z}\\
=& \frac{(1+z)P^{\bar{R}}_{\bar{R}/\bar{I}}-1}{z}
= \frac{(1+z)(1+zP_{\bar{I}}^{\bar{R}}(z))-1}{z}\\
=& \frac{(1+z)(1+h_{2}zP_{\mathbb{K}}^{\bar{R}}(z))-1}{z}\\
=& \frac{(1+z)(1+h_{2}zP_{\mathbb{K}}^{R}(z)/(1+z))-1}{z}\\
=& \frac{h_{2}zP^{R}_{\mathbb{K}}(z)+z}{z}=1+h_{2}P_{\mathbb{K}}^{R}(z).
\end{align*}
Also, when $\Gamma_{2}=\mathbb{N}$, by Lemma \ref{generator} and Lemma \ref{iso k}, 
$\bar{I}\cong \mathbb{K}$ implies $P^{R}_{\bar{D}_{1}}(z)=1+P_{\mathbb{K}}^{R}(z)$.
\medskip

\item When $\Gamma_{1} \neq \mathbb{N}$, we have 
$\bar{J}=(\bar{v}^{\delta_{1}+n_{e}-1}\bar{u}^{\gamma_{1}},\dots,\bar{v}^{\delta_{h_{1}}+n_{e}-1}\bar{u}^{\gamma_{h_{1}}})\cong \mathbb{K} \oplus \mathbb{K} \dots\oplus \mathbb{K} \, 
(h_{1}\textrm{-copies})$. By the same argument as in (i), 
we have $P_{\bar{D}_{2}}^{R}(z)=1+h_{1}P_{\mathbb{K}}^{R}(z)$. In the case when 
$\Gamma_{1}=\mathbb{N}$, by Lemma \ref{generator} and Lemma \ref{iso k}, 
$\bar{J}\cong \mathbb{K}$ implies $P^{R}_{\bar{D}_{2}}(z)=1+P_{\mathbb{K}}^{R}(z)$.
\medskip

\item By Lemma \ref{generator} and Lemma \ref{ideal}, $\mathrm{Der}_{\mathbb{K}}(R)$ 
is isomorphic to the ideal $I_{\mathrm{Der}_{\mathbb{K}}(R)}$ in 
$\mathbb{K}[\overline{\Gamma}]$ and 
\[
   I_{\mathrm{Der}_{k}(R)}= 
\begin{cases}
    (u^{n_{e}},v^{\beta_{1}}u^{\alpha_{1}-1+n_{e}},\dots,v^{\beta_{h_{2}}}u^{\alpha_{h_{2}}-1+n_{e}},v^{n_{e}}, v^{\delta_{1}+n_{e}-1}u^{\gamma_{1}},\dots,v^{\delta_{h_{1}}+n_{e}-1}u^{\gamma_{h_{1}}}) & \text{if } \Gamma_{1} \neq \mathbb{N}, \Gamma_{2} \neq \mathbb{N},
\\
    (u^{n_{e}},v^{\beta_{1}}u^{\alpha_{1}-1+n_{e}},\dots,v^{\beta_{h_{2}}}u^{\alpha_{h_{2}}-1+n_{e}},v^{n_{e}},v^{n_{e}-1}u^{1+e'n_{e}}) & \text{if } \Gamma_{1} = \mathbb{N}, \Gamma_{2} \neq \mathbb{N},\\
(u^{n_{e}},u^{n_{e}-1}v^{1+c'n_{e}},v^{n_{e}},v^{\delta_{1}+n_{e}-1}u^{\gamma_{1}},\dots,v^{\delta_{h_{1}}+n_{e}-1}u^{\gamma_{h_{1}}}) & \text{if } \Gamma_{1} \neq \mathbb{N}, \Gamma_{2} = \mathbb{N}, \\
(u^{n_{e}},u^{n_{e}-1}v^{1+c'n_{e}},v^{n_{e}},v^{n_{e}-1}u^{1+e'n_{e}}) & \text{otherwise.}
\end{cases}
\]
Clearly the image of $I_{\mathrm{Der}_{k}(R)}$ in $\bar{R}$ is
\[
= 
\begin{cases}
((\bar{v}^{\beta_{1}}\bar{u}^{\alpha_{1}-1+n_{e}},\dots,\bar{v}^{\beta_{h_{2}}}\bar{u}^{\alpha_{h_{2}}-1+n_{e}},\bar{v}^{\delta_{1}+n_{e}-1}\bar{u}^{\gamma_{1}},\dots,\bar{v}^{\delta_{h_{1}}+n_{e}-1}\bar{u}^{\gamma_{h_{1}}})\\
 \cong \mathbb{K} \oplus \mathbb{K} \oplus \cdots \oplus \mathbb{K} \, ((h_{1}+h_{2}) \, \rm{copies}), & \text{if } \Gamma_{1} \neq \mathbb{N}, \Gamma_{2} \neq \mathbb{N},\\
(\bar{v}^{\beta_{1}}\bar{u}^{\alpha_{1}-1+n_{e}},\dots,\bar{v}^{\beta_{h_{2}}}\bar{u}^{\alpha_{h_{2}}-1+n_{e}},\bar{v}^{n_{e}-1}\bar{u}^{1+e'n_{e}})    \\
 \cong \mathbb{K} \oplus \mathbb{K} \oplus \cdots \oplus \mathbb{K} \, 
 ((h_{2}+1) \, \rm{copies})         & \text{if } \Gamma_{1} = \mathbb{N}, \Gamma_{2} \neq \mathbb{N},\\
(\bar{u}^{n_{e}-1}\bar{v}^{1+c'n_{e}},\bar{v}^{\delta_{1}+n_{e}-1}\bar{u}^{\gamma_{1}},\dots,\bar{v}^{\delta_{h_{1}}+n_{e}-1}\bar{u}^{\gamma_{h_{1}}})\\
\cong \mathbb{K} \oplus \mathbb{K} \oplus \cdots \oplus \mathbb{K} \, 
((1+h_{1}) \, \rm{copies}) & \text{if } \Gamma_{1} \neq \mathbb{N}, \Gamma_{2} = \mathbb{N},\\
(\bar{u}^{n_{e}-1}\bar{v}^{1+c'n_{e}},\bar{v}^{n_{e}-1}\bar{u}^{1+e'n_{e}}) \\
\cong \mathbb{K} \oplus \mathbb{K}, & \text{otherwise}.
\end{cases}
\]
Now, using Lemma \ref{COM P} with $M=R/I_{\mathrm{Der}_{\mathbb{K}}(R)}$ and 
the argument of (i), we get 
$$P_{\mathrm{Der}_{\mathbb{K}}(R)}^{R}(z)=P_{I_{\mathrm{Der}_{\mathbb{K}}(R)}}^{R}(z)=1+(h_{1}+h_{2})P_{\mathbb{K}}^{R}(z)$$(when $\Gamma_{i}=\mathbb{N}$, so that $\vert PF(\Gamma_{i})\vert =h_{i}=1$).\qed
\end{enumerate}

\begin{corollary}
 Let $R=\mathbb{K}[u^{n_{e}},v^{n_{1}}u^{n_{e}-n_{1}},\dots,v^{n_{e-1}}u^{n_{e}-n_{e-1}},v^{n_{e}}] $ be a Cohen-Macaulay ring. Suppose $\Gamma_{1}=\langle n_{1},\dots,n_{e} \rangle $ be a numerical semigroup with $\vert \mathrm{PF}(\Gamma_{1}) \vert=h_{1}$ and 
 $\Gamma_{2}=\langle n_{e}-n_{1},\dots,n_{e}-n_{e-1},n_{e} \rangle$ be the 
 numerical semigroup with $\vert \mathrm{PF}(\Gamma_{2}) \vert=h_{2}$. 
 Then the Betti numbers of $\mathrm{Der}_{\mathbb{K}}(R)$ over $R$ are 
 $(1+h_{1}h_{2})$-times the Betti numbers of $\mathbb{K}$ over $R$.
\end{corollary} 

\proof Let $P_{\mathbb{K}}^{R}(z)= \displaystyle\sum _{i \geq 0}\beta_{i}^{R}(\mathbb{K})z^{i}$ 
be the Poincar\'{e} series of $\mathbb{K}$ over $R$. From Lemma \ref{Poincare-1}, 
we have 
$$P_{\mathrm{Der}_{\mathbb{K}}(R)}^{R}(z) = \displaystyle\sum _{i \geq 0}\beta_{i}^{R}(\mathrm{Der}_{\mathbb{K}}(R))z^{i} = 1+(h_{1}+h_{2})P_{\mathbb{K}}^{R}(z) = 1+(h_{1}+h_{2})\displaystyle\sum _{i \geq 0}\beta_{i}^{R}(\mathbb{K})z^{i}.$$
Comparing the coefficients of the above expression,  
$\beta_{i}^{R}(\mathrm{Der}_{\mathbb{K}}(R))=1+(h_{1}+h_{2})\beta_{i}^{R}(\mathbb{K})$. \qed 

\section{Derivation Module of Projective Closure of the Arslan Curve}
In this section we find the explicit description of the module of derivations of 
the Arslan curve. Let us first recall Arslan's example of monomial curves 
in the affine space $\mathbb{A}^{4}$. Arslan defined the numerical semigroups $\Gamma_{A_{h}}=
\Gamma( h(h+1),h(h+1)+1,(h+1)^{2},(h+1)^{2}+1)$, for $h\geq 2$. Assume that minimal 
generator of $\Gamma_{A_{h}}$ are denoted by $a_{1h},a_{2h},a_{3h},a_{4h}$ respectively.

\begin{lemma}[Proposition 3.2 \cite{Arslan}] 
Let us consider the following polynomials, for $h\geq 2$ : 
\medskip

\noindent $ w=x_{1}x_{4}-x_{2}x_{3}$;\\[2mm]
$ g_{i}=x_{1}^{h-i}x_{3}^{i+1}-x_{2}^{h-i+1}x_{4}^{i},$ \, for \, $ 0 \leq i \leq h$;\\[2mm]
$ q_{j}=x_{3}^{j}x_{4}^{h-j}-x_{1}^{j+1}x_{2}^{h-j},$ for $ 0 \leq j \leq h $.
\medskip

\noindent The set $G_{h}=\{ w,g_{i},q_{j}\mid 1\leq i,j\leq h\}$ 
is a minimal genrating set of the defining ideal of the Arslan curves $\Gamma_{A_{h}}$.
\end{lemma}
\medskip

Let $\Gamma=\langle n_{1},\dots,n_{e}\rangle$ be a numerical semigroup. Given $0 \neq s \in \Gamma,$ the set of lengths of $s$ in $\Gamma$ is defined as $\mathcal{T}(s)=\Bigg\{\sum_{i=1}^{e} r_{i} \mid s=\sum_{i=1}^{e}r_{i}n_{i},r_{i}\geq 0 \Bigg\}$.

\begin{definition}(\cite{Homogeneous})
A subset $T \subset \Gamma$ is called \emph{homogeneous} if either it is empty or $\mathcal{T}(s)$ is singleton for all $0 \neq s \in T$. The numerical semigroup $\Gamma$ is called homogeneous, when the Ap\'{e}ry set $\mathrm{AP}(\Gamma, n_{1})$ is homogeneous.
\end{definition}

\begin{theorem}
The numerical semigroup $\Gamma_{A_{h}}$ is homogeneous.
\end{theorem}

\proof Since $x_{1}$ belongs to the support of a non-homogeneous element $q_{j}\in G_{h}$, 
by Corollary 3.10 in \cite{Homogeneous} we have $\mathrm{Ap}(\Gamma_{A_{h}},a_{1h})$ is homogeneous 
and hence $\Gamma_{A_{h}}$ is homogeneous.

\begin{lemma}
Let $h \geq 2$. Let $\mathrm{Ap}(\Gamma_{A_{h}},a_{1h})$ denote the Ap\'{e}ry set of $\Gamma_{A_{h}}$, 
with respect to the element $a_{1h}$. Then 
\begin{align*}
\mathrm{Ap}(\Gamma_{A_{h}},a_{1h})=&  \{ia_{2h} \mid 0 \leq i \leq h\}\\
&\cup \{ja_{3h} \mid 1 \leq j \leq h-1\}\\
&\cup \{la_{4h} \mid 1 \leq l \leq h-1\} \\
&\cup \{ \gamma a_{2h}+\nu a_{4h} \mid 2 \leq \gamma+\nu \leq  h, \gamma \neq 0, \nu \neq 0\} \\
&\cup \{\alpha a_{3h}+\beta a_{4h} \mid 2 \leq \alpha+\beta \leq  h-1, \alpha \neq 0, \beta \neq 0\}
\end{align*}
\end{lemma}
\proof Follows from Example 3.13 in \cite{Type-Arslan}. \qed

\begin{definition}
Let $\Gamma$ be a numerical semigroup; $x \in \mathbb{Z}$ is a pseudo-Frobenius number if 
$\, x \notin \Gamma$ and $ x+s \in \Gamma \,\,\ \text{for all}\,\,s \in \Gamma \setminus \{0\}$. 
The set of pseudo-Frobenius numbers is denoted by $\mathrm{PF}(\Gamma)$.
\end{definition}

\begin{lemma}
Let $h \geq 2$. Then
$$\mathrm{PF}(\Gamma_{A_{h}})= \{(h-i)a_{3h}+(i-1)a_{4h}-a_{1h} \mid 1 \leq i \leq h-1\} \cup \{(h-j)a_{2h}+ja_{4h}-a_{1h} \mid 0 \leq j \leq h-1\}.$$
\end{lemma} 
\proof See Example 3.13 in \cite{Type-Arslan}. \qed
\medskip

\begin{theorem}\label{Der-Arslan}
Let $R=\mathbb{K}[u^{a_{4h}},v^{a_{1h}}u^{a_{4h}-a_{1h}},v^{a_{2h}}u^{a_{4h}-a_{2h}},v^{a_{3h}}u^{a_{4h}-a_{3h}},v^{a_{4h}}]$ be the projective closure of the Arslan curve $\Gamma_{A_{h}}$. 
For $1 \leq i \leq h-1$, $ 0 \leq j \leq h-2$, the derivation module $\mathrm{Der}_{\mathbb{K}}(R)$ 
is generated by $\displaystyle\bigcup_{i=1}^{6}A_{i}$, where
\begin{itemize}
\item $A_{1}=\{v^{1+ha_{4h}}\frac{\partial} {\partial v}\}$;
\item  $A_{2}=\{v^{(h-i)a_{3h}+(i-1)a_{4h}-a_{1h}+1}u^{h(h+1)+h-i}\frac{\partial} {\partial u}\}$;
\item  $A_{3}= \{v^{(h-j)a_{2h}+ja_{4h}-a_{1h}+1}u^{(h+1)(h-j-1)-1}\frac{\partial} {\partial u}\}$;
\item  $A_{4}=\{v^{a_{2h}+(h-1)a_{4h}-a_{1h}+1}u^{(h+1)^{2}}\frac{\partial} {\partial u}\}$;
\item  $A_{5}=\{u \frac{\partial} {\partial u}\}$;
\item  $A_{6}=\{v \frac{\partial} {\partial v}\}$.
\end{itemize}
Hence, $\mu(\mathrm{Der}_{\mathbb{K}}(R))=2h+2$.
\end{theorem}

\proof We use Theorem \ref{generator}. We have $\Gamma_{A_{h}}=\langle h(h+1),h(h+1)+1,(h+1)^{2},(h+1)^{2}+1 \rangle$ and $\Gamma_{2}= \mathbb{N}$. By Proposition 5.2, \cite{CMC}, 
$R$ is Cohen-Macaulay, therefore we have $\overline{\Gamma}_{A_{h}}=(\Gamma_{A_{h}} \times \Gamma_{2})\cap L$ 
(See \cite{Niesi}, Theorem 4.6). 

We claim that there does not exist $i>0$ such that $1+(h-i)((1+h)^{2}+1) \in \Gamma_{A_{h}}$.
We prove our claim by induction. For $i=1$, we have $1+(h-1)((1+h)^{2}+1)= (h+1)(h(h+1))-1-(h(h+1)) \notin \Gamma_{1}$, since  $(h+1)(h(h+1))-1 \in \mathrm{Ap}(\Gamma_{A_{h}},a_{1h})$. 
We now assume $1+(h-i)((1+h)^{2}+1) \notin \Gamma_{A_{h}}$. For $i+1$, suppose $1+(h-(i+1))((1+h)^{2}+1)=1+(h-i)((1+h)^{2}+1)-((1+h)^{2}+1) \in \Gamma_{A_{h}}$. Since $((1+h)^{2}+1)\in \Gamma_{A_{h}}$, 
hence $1+(h-i)((1+h)^{2}+1) \in \Gamma_{A_{h}}$, which is a contradiction. Therefore, $h$ 
is the least natural number such that $1+h((1+h)^{2}+1) \in \Gamma_{A_{h}}$ and we get the following:
\begin{enumerate}
\item[(a)] $(1+h((h+1)^{2}+1),-1)+(0,a_{4h}) \in \overline{\Gamma}_{A_{h}}$.
\medskip

\item[(b)] $(1+h((h+1)^{2}+1),-1)+(a_{1h},a_{4h}-a_{1h}) = (h((h+1)^{2}+1)+h(h+1)+1,a_{4h}-a_{1h}-1) \in (\Gamma_{1}\times \Gamma_{2})\cap L= \overline{\Gamma}_{A_{h}}$. 
\medskip

\item[(c)] $1+h((h+1)^{2}+1),-1)+(a_{2h},a_{4h}-a_{2h}) = (h^{2}(h+1)+h(h+1),a_{4h}-a_{2h}-1)\in (\Gamma_{1}\times \Gamma_{2})\cap L= \overline{\Gamma}_{A_{h}}$. 
\medskip

\item[(d)] $1+h((h+1)^{2}+1),-1)+(a_{3h},a_{4h}-a_{3h}) = 
(1+h((1+h)^{2}+1)+(1+h)^{2},a_{4h}-a_{3h}-1)\in (\Gamma_{1}\times \Gamma_{2})\cap L= \overline{\Gamma}_{A_{h}}$. 
\end{enumerate}
Thus, $h$ is the least natural number such that $(1+ha_{4h},-1)+(n,a_{4h}-n)\in \overline{\Gamma}_{A_{h}}$, 
for each $n \in \{0,a_{1h},a_{2h},a_{3h}\}$. Therefore, 
we have $D_{1}=\{v^{1+ha_{4h}}\frac{\partial} {\partial v}\}$.
\medskip

Let $t_{i}=(h-i)a_{3h}+(i-1)a_{4h}-a_{1h}$, where $1 \leq i \leq h-1$.
Now we find the least natural number $\gamma_{t_{i}}$ such that $(t_{i},\gamma_{t_{i}}) \in (\mathrm{PF}(\Gamma_{A_{h}})\times \Gamma_{2})\cap L$. We claim that $\gamma_{t_{i}}=h(h+1)+h-i$. We consider
 \begin{align*}
&t_{i}-(h-2)((h+1)^{2}+1)\\
=&(h-i)(h+1)^{2}+(i-1)((h+1)^{2}+1)-h(h+1)-(h-2)((h+1)^{2}+1)\\
=&(h+1)[(h-i)(h+1)+(i-h+1)(h+1)-h]+i-h+1\\
=&i+2 > 0.
\end{align*}
Therefore, $ t_{i}>(h-2)((h+1)^{2}+1)$ and $t_{i}+ h(h+1)+h-i= (h-1)((h+1)^{2}+1)$. Hence, 
$\gamma_{t_{i}}=h(h+1)+h-i$ and it is clear that 
$(t_{i},\gamma_{t_{i}}) \in (PF(\Gamma_{A_{h}})\times \Gamma_{2})$. Thus, 
$(t_{i},\gamma_{t_{i}}) \in (PF(\Gamma_{A_{h}})\times \Gamma_{2})\cap L$. Since $t_{i}\in \mathrm{PF}(\Gamma_{A_{h}})$, 
we get $t_{i}+n \in \Gamma_{A_{h}}$ and $(t_{i},\gamma_{t_{i}})+(n,a_{4h}-n)\in \Gamma_{A_{h}}$, 
for each $n \in \{a_{1h},a_{2h},a_{3h},a_{4h}\}$.
\medskip

Let $y_{j}=(h-j)a_{2h}+ja_{4h}-a_{1h}$, for $0 \leq j \leq h-2$. We find the least natural number $\gamma_{y_{j}}$ such that $(y_{j},\gamma_{y_{j}}) \in (PF(\Gamma_{A_{h}})\times \Gamma_{2})\cap L$. It is obvious that $(y_{j},\gamma_{y_{j}}) \in (PF(\Gamma_{A_{h}})\times \Gamma_{2})$. We claim that $\gamma_{y_{j}}=h(h+1)-(j+1)(h+1)-1$. 
We observe that
\begin{align*}
&y_{j}-(h-2)((h+1)^{2}+1)\\
=&(h-j)(h(h+1)+1)+j((h+1)^{2}+1)-h(h+1)-(h-2)((h+1)^{2}+1)\\
=&(h+1)[h(h-j)+j(h+1)-h-(h-2)(h+1)]+h-j+j-h+2\\
=&(h+1)(j+2)+2>0.
\end{align*}
Therefore $ y_{j}>(h-2)((h+1)^{2}+1)$ and $y_{j}+h(h+1)-(j+1)(h+1)-1=(h-1)((h+1)^{2}+1)$, 
for $0 \leq j \leq h-2$. Hence, $\gamma_{y_{j}}=h(h+1)-(j+1)(h+1)-1$ and $(y_{j},\gamma_{y_{j}})\in (PF(\Gamma_{A_{h}})\times \Gamma_{2})$. Thus, 
$(y_{j},\gamma_{y_{j}}) \in (PF(\Gamma_{A_{h}})\times \Gamma_{2})\cap L$. Since $y_{j}\in \mathrm{PF}(\Gamma_{A_{h}})$, we get $y_{j}+n \in \Gamma_{A_{h}}$ and 
$(y_{j},\gamma_{y_{j}})+(n,a_{4h}-n)\in \Gamma_{A_{h}}$, for each $n \in \{a_{1h},a_{2h},a_{3h},a_{4h}\}$.
\medskip

Next we want to find the least natural number $\gamma_{w_{h-1}}$, 
such that $(a_{2h}+(h-1)a_{4h}-a_{1h},\gamma_{w_{h-1}})\in (PF(\Gamma_{A_{h}})\times \Gamma_{2})$. 
We note that 
\begin{align*}
 w_{h-1}+(h+1)^{2}&=(h(h+1)+1)+(h-1)((h+1)^{2}+1)-h(h+1)+(h+1)^{2}\\
 &=(h-1)((h+1)^{2}+1). 
\end{align*}
Therefore, $w_{h-1}>(h-2)((h+1)^{2}+1)$. Hence, $\gamma_{w_{h-1}}=(h+1)^{2}$ and $(w_{h-1},\gamma_{w_{h-1}})\in (PF(\Gamma_{A_{h}})\times \Gamma_{2})$. Thus, $(w_{h-1},\gamma_{w_{h-1}}) \in (PF(\Gamma_{A_{h}})\times \Gamma_{2})\cap L$. Since $w_{h-1}\in \mathrm{PF}(\Gamma_{A_{h}})$, we get $w_{h-1}+n \in \Gamma_{A_{h}}$ and $(w_{h-1},\gamma_{w_{h-1}})+(n,a_{4h}-n)\in \Gamma_{A_{h}}$, for each $n \in \{a_{1h},a_{2h},a_{3h},a_{4h}\}$. 
\medskip

Therefore, using Theorem \ref{generator}, we have for $1 \leq i \leq h-1$, $ 0 \leq j \leq h-2$,
\begin{align*}
D_{2}&=  \{v^{(h-i)a_{3h}+(i-1)a_{4h}-a_{1h}+1}u^{h(h+1)+h-i}\frac{\partial} {\partial u}\}\\
         &\cup \{v^{(h-j)a_{2h}+ja_{4h}-a_{1h}+1}u^{(h+1)(h-j-1)-1}\frac{\partial} {\partial u}\}\\
         &\cup \{v^{a_{2h}+(h-1)a_{4h}-a_{1h}+1}u^{(h+1)^{2}}\}. \qed
\end{align*}

\begin{remark}
Note that $| \mathrm{PF}(\mathrm{pr}_{1}(\overline{\Gamma}_{A_{h}}))|=| \mathrm{PF}(\Gamma_{A_{h}})|=2h-1 $ and $| \mathrm{PF}(\mathrm{pr}_{2}(\overline{\Gamma}_{A_{h}}))|=1$. Therefore, by Theorem \ref{Poincare-1}, $P_{\mathrm{Der}_{\mathbb{K}}(R)}^{R}(z)= 1+(2h)P_{\mathbb{K}}^{R}(z)$.
\end{remark}

\section{Derivation module of the Backelin Curves}
Backelin defined the numerical semigroup $\Gamma_{B_{nr}}=\Gamma( s,s+3,s+3n+1,s+3n+2 )$, 
for $n \geq 2, r\geq 3n+2 $ and $ s=r(3n+2)+3$ in $\mathbb{A}_{\mathbb{K}}^{4}$. For the sake 
of convenience, let us rename 
the minimal generators of the numerical semigroup $\Gamma_{B_{nr}}$ as follows: 
$s = m_{1_{nr}}$, $s+3 = m_{2_{nr}}$, $s+3n+1 = m_{3_{nr}}$, and $s+3n+2 = m_{4_{nr}}$. 
Let the defining ideal be denoted by $\mathfrak{p}(m_{1_{nr}},m_{2_{nr}},m_{3_{nr}},m_{4_{nr}})$. 

\begin{lemma}[Theorem 2.2 \cite{Backelin}] \label{generator-Backelin}
The defining ideal $\mathfrak{p}(m_{1_{nr}},m_{2_{nr}},m_{3_{nr}},m_{4_{nr}})$ of monomial curve associated to $\Gamma_{B_{nr}}$ is minimally generated by following polynomials
\begin{enumerate}[(i)]
\item $f_{1}=x_{2}x_{3}^{3}-x_{1}x_{4}^{3}$
\item $f_{(2,i)}=x_{1}^{n-i}x_{3}^{3i-1}-x_{2}^{n-i+1}x_{4}^{3i-2},\,\,1\leq i \leq n$
\item $f_{(3,j)}=x_{1}^{r-n+3+j}x_{2}^{n-1-j}-x_{3}^{2+3j}x_{4}^{r-1-3j},\,\,0 \leq j \leq n-1$
\item $f_{(4,j)}=x_{1}^{r-2n+3+j}x_{2}^{2n-j}-x_{3}^{3j+1}x_{4}^{r+1-3j}, \,\, 0 \leq j \leq n-1$
\item $f_{5}=x_{1}^{r-n+2}x_{2}^{n}x_{3}-x_{4}^{r+2}$
\item $f_{6}=x_{2}^{n+1}x_{3}-x_{1}^{n}x_{4}^{2}$
\item $f_{7}=x_{2}^{2n+1}-x_{1}^{2n-1}x_{3}x_{4}$.
\end{enumerate}
\end{lemma}

\begin{theorem}
The numerical semigroup $\Gamma_{B_{nr}}$ is homogeneous.
\end{theorem}

\proof From Lemma \ref{generator-Backelin}, it is easy to see that $x_{1}$ belongs to the support of the non-homogeneous generator of $\mathfrak{p}(m_{1_{nr}},m_{2_{nr}},m_{3_{nr}},m_{4_{nr}})$. Hence, 
$\Gamma_{B_{nr}}$ is homogeneous by \cite{Homogeneous} .

\begin{lemma}[Proposition 1.8 \cite{Pseduo-Backelin}]
Let $n\geq 2$. The pseudo-Frobenius set of $\Gamma_{B_{nr}}$ is given by 
$\mathrm{PF}(\Gamma_{B_{nr}}) = F_{1} \cup F_{2} \cup F_{3} \cup F_{4}$, 
where 
\begin{align*}
F_{1} := & \{(n-l)m_{1_{nr}}+(3l-2)m_{3_{nr}}-m_{4_{nr}} \mid 2\leq l \leq n\};\\
F_{2} := & \{(r-(n+m)+3)m_{1_{nr}}+(n+m-1)m_{2_{nr}}-m_{4_{nr}} \mid 1\leq m \leq n\};\\
F_{3} := & \{(r-m+2)m_{1_{nr}}+(m-1)m_{2_{nr}}+m_{3_{nr}}-m_{4_{nr}} \mid 1\leq m \leq n\};\\
F_{4} := & \{(r-n+1)m_{1_{nr}}+n\cdot m_{2_{nr}}+m_{3_{nr}}-m_{4_{nr}}\} \cup \\
& \{(n-2)m_{1_{nr}}+n\cdot m_{2_{nr}}+2m_{3_{nr}}-m_{4_{nr}}\} \cup 
\{(r-2n+2)m_{1_{nr}}+2n\cdot m_{2_{nr}}-m_{4_{nr}}\}. 
\end{align*}
\end{lemma}

\begin{theorem}\label{Der-Backelin}
Let $R=\mathbb{K}[u^{m_{4_{nr}}},v^{m_{1_{nr}}}u^{m_{4_{nr}}-m_{1_{nr}}},v^{m_{2_{nr}}}u^{m_{4_{nr}}-m_{2_{nr}}},v^{m_{3_{nr}}}u^{m_{4_{nr}}-m_{3_{nr}}},v^{m_{4_{nr}}}]$ be the 
projective closure of the Backelin curve $\Gamma_{B_{nr}}$. Then $\mathrm{Der}_{\mathbb{K}}(R)$ 
is generated by $\displaystyle\bigcup_{i=1}^{8}B_{i}$, where
\begin{itemize}
\item $B_{1}=\{v^{1+(r+1)m_{4_{nr}}}\frac{\partial}{\partial u}\}$;
\item $B_{2}= \{v^{(n-l)m_{1_{nr}}+(3l-2)m_{3_{nr}}-m_{4_{nr}}+1}u^{3n^{2}+2n-3ln+l-2}\frac{\partial}{\partial v}\mid l=2,\dots,n,\}$;
\item $B_{3}= \{v^{(r-(n+m)+3)m_{1_{nr}}+(n+m-1)m_{2_{nr}}-m_{4_{nr}}+1}u^{s+3n+1-3(m-1)}\frac{\partial}{\partial v}\mid m=1,\dots,n\}$;
\item $B_{4}= \{v^{(r-m+2)m_{1_{nr}}+(m-1)m_{2_{nr}}+m_{3_{nr}}-m_{4_{nr}}+1}u^{s+3n-3(m-1)}\frac{\partial}{\partial v}\mid m=1,\dots,n\}$;
\item $B_{5}= \{v^{(r-n+1)m_{1_{nr}}+n\cdot m_{2_{nr}}+m_{3_{nr}}-m_{4_{nr}}+1}u^{s}\frac{\partial}{\partial v}\}$;
\item $B_{6}=\{v^{(n-2)m_{1_{nr}}+n\cdot m_{2_{nr}}+2m_{3_{nr}}-m_{4_{nr}}+1}u^{6n^{2}-5n-2}\frac{\partial}{\partial v}\} $;
\item $B_{7}=\{v^{(r-2n+2)m_{1_{nr}}+2n\cdot m_{2_{nr}}-m_{4_{nr}}+1}u^{s+1}\frac{\partial}{\partial v}\} $;
\item $B_{8}= \{u\frac{\partial}{\partial u},v\frac{\partial}{\partial v}\} $.
\end{itemize}
 Hence, $\mu(\mathrm{Der}_{\mathbb{K}}(R))=3n+5$.
\end{theorem}

\proof  We have $\Gamma_{B_{nr}}=\langle s,s+3,s+3n+1,s+3n+2 \rangle$ and $\Gamma_{2}=\mathbb{N}$. By 
Theorem 2.8 \cite{Backelin}, $R$ is Cohen-Macaulay. Therefore, it 
follows from Theorem 4.6 \cite{Niesi} that $\overline{\Gamma}_{B_{nr}}=(\Gamma_{B_{nr}}\times \Gamma_{2})\cap L$.
\medskip

We claim that there does not exist any $i>0$, such that $1+(r+1-i)(s+3n+2) \in \Gamma_{B_{nr}}$. 
We prove our claim by induction. For $i=1$, we have $1+r(s+3n+2)=1+(r+1)(s+3n+2)-(s+3n+2) \notin \Gamma_{B_{nr}}$ as $1+(r+1)(s+3n+2)\in \mathrm{Ap}(\Gamma_{B_{nr}},m_{4_{nr}})$. 
\smallskip

By induction hypothesis we assume that,  $1+(r+1-i)(s+3n+2) \notin \Gamma_{B_{nr}}$. Now for $i+1$, suppose $1+(r+1-(i+1))(s+3n+2) \in \Gamma_{B_{nr}}$. 
But we have  $$1+(r+1-i)(s+3n+2)=1+(r+1-(i+1))(s+3n+2)+(s+3n+2) \in \Gamma_{B_{nr}}$$ and which gives a contradiction of the induction hypothesis. 
Hence there does not exist any $i>0$ for which $1+(r+1-i)(s+3n+2) \in \Gamma_{B_{nr}}$. Thus, $r+1$ is the least natural number such that $(1+(r+1)m_{4_{nr}},-1)+(n,m_{4_{nr}}-n)\in \overline{\Gamma}_{B}$ for each $n \in \{0,m_{1_{nr}},m_{2_{nr}},m_{3_{nr}}\}$. 
Thus, $r+1$ is the least natural number such that $1+(r+1)(s+3n+2) \in \Gamma_{B_{nr}}$. 
It is easy to check that, 
\begin{itemize}
\item $(1+(r+1)(s+3n+2),-1)+(0,m_{4_{nr}})\in \overline{\Gamma}_{B_{nr}}$;
\item $(1+(r+1)(s+3n+2),-1)+(m_{1_{nr}},m_{4_{nr}}-m_{1_{nr}})\in \overline{\Gamma}_{B_{nr}}$;
\item $(1+(r+1)(s+3n+2),-1)+(m_{2_{nr}},m_{4_{nr}}-m_{2_{nr}})\in \overline{\Gamma}_{B_{nr}}$;
\item $(1+(r+1)(s+3n+2),-1)+(m_{3_{nr}},m_{4_{nr}}-m_{3_{nr}})\in \overline{\Gamma}_{B_{nr}}$.
\end{itemize}
Thus, $r+1$ is the least natural number such that 
$(1+(r+1)m_{4_{nr}},-1)+(n,m_{4_{nr}}-n)\in \overline{\Gamma}_{B_{nr}}$, 
for each $n \in \{0,m_{1_{nr}},m_{2_{nr}},m_{3_{nr}}\}$. Therefore, according to the 
notations of Theorem \ref{generator}, $D_{1}=\{v^{1+(r+1)m_{4_{nr}}}\frac{\partial} {\partial v}\}$.

Let $p_{l}=(n-l)m_{1_{nr}}+(3l-2)m_{3_{nr}}-m_{4_{nr}} \in \mathrm{PF}(\Gamma_{B_{nr}})$, 
for $2  \leq l \leq n$. We now find the least natural number $\gamma_{p_{l}}$, 
such that $(p_{l},\gamma_{p_{l}})\in (\mathrm{PF}(\Gamma_{B_{nr}})\times \Gamma_{2})\cap L$. 
We claim that, $\gamma_{p_{l}}=3n^{2}+2n-3ln+l-2$. It is easy to see that
{\scriptsize
\begin{align*}
p_{l}-(n+2l-4)m_{4_{nr}} & = (n-l)s+(3l-2)(s+3n+1)-(s+3n+2)-(n+2l-4)(s+3n+2)\\
&=s[n-l+3l-2-1-n-2l-4]+(3l-2)(3n+1)-3n-2-(2l-4)(3n+2)\\
&=s+3nl+3n-l+8>0.
\end{align*}
}
Therefore, $ p_{l}>(n+2l-4)m_{4_{nr}}$, for $l=2,\ldots,n$, \, and  \, $p_{l}+3n^{2}+2n-3ln+l-2=(n+2l-3)m_{4_{nr}}$. Hence, $\gamma_{p_{l}}=3n^{2}+2n-3ln+l-2$ and therefore $(p_{l},\gamma_{p_{l}})=(\mathrm{PF}(\Gamma_{B_{nr}})\times \Gamma_{2})\cap L$.

Let $q_{m}=(r-(n+m)+3)m_{1_{nr}}+(n+m-1)m_{2_{nr}}-m_{4_{nr}} \in \mathrm{PF}(\Gamma_{B_{nr}})$, 
where $1\leq m \leq n$. We consider 
{\scriptsize
\begin{align*}
q_{m}-r\cdot m_{4_{nr}} & = (r-(n+m)+3)m_{1_{nr}}+(n+m-1)m_{2_{nr}}-m_{4_{nr}}-r\cdot m_{4_{nr}}\\
=& (r-(n+m)+3)s+(n+m-1)(s+3)-(r+1)(s+3n+2)\\
=&s[r-(n+m)+3+n+m-1-(r+1)]+3(n+m)-3+(r+1)(3n+2)\\
=& s+3(n+m)+3n(r+1)+3n-1>0.
\end{align*}
}
Therefore, $q_{m}>r\cdot m_{4_{nr}}$, where $1 \leq m \leq n$, \, and \, $q_{m}+(s+3n+1-3(m-1))=(r+1)m_{4_{nr}}$. Hence, $\gamma_{q_{m}}=s+3n+1-3(m-1)$ and $(q_{m},\gamma_{q_{m}})\in (\mathrm{PF}(\Gamma_{B_{nr}})\times \Gamma_{2}) \cap L$. Since $q_{m} \in \mathrm{PF}(\Gamma_{B_{nr}})$, we get $q_{m}+n \in \Gamma_{B_{nr}}$ and $(q_{m},\gamma_{q_{m}})+(n,m_{4_{nr}}-n)\in \overline{\Gamma}_{B}$, for each $\{m_{1_{nr}},m_{2_{nr}},m_{3_{nr}},m_{4_{nr}}\}$.

Suppose $z_{m}=(r-m+2)m_{1_{nr}}+(m-1)m_{2_{nr}}+m_{3_{nr}}-m_{4_{nr}}\in \mathrm{PF}(\Gamma_{B_{nr}})$, 
where $1 \leq m \leq n$. Consider 
{\scriptsize
\begin{align*}
z_{m}-r\cdot m_{4_{nr}} & =  (r-m+2)s+(m-1)(s+3)+(s+3n+1)-(s+3n+2)-r\cdot (s+3n+2)\\
=& s(r-m+2+m-1-r)-1+3(m-1)-r(3n+2)\\
=& s+3m-4-r(3n+2)=3m-1>0.
\end{align*}
}
Therefore, $ z_{m}>r\cdot m_{4_{nr}}$, where $1 \leq m \leq n$, \, and 
$Z_{m}+(s+3n-3(m-1)=(r+1)m_{4_{nr}}$. Thus, we get $\gamma_{z_{m}}=s+3n-3(m-1)$ 
and it is clear that $(z_{m},\gamma_{z_{m}})\in (\mathrm{PF}(\Gamma_{B_{nr}})\times \Gamma_{2}) \cap L$. 
Since $Z_{m} \in \mathrm{PF}(\Gamma_{B_{nr}})$, we get $z_{m}+n \in \Gamma_{B_{nr}}$ and $(z_{m},\gamma_{z_{m}})+(n,m_{4_{nr}}-n)\in \overline{\Gamma}_{B}$, for each $\{m_{1_{nr}},m_{2_{nr}},m_{3_{nr}},m_{4_{nr}}\}$.

Now we take $P=(r-n+1)m_{1_{nr}}+n\cdot m_{2_{nr}}+m_{3_{nr}}-m_{4_{nr}} \in \mathrm{PF}(\Gamma_{B_{nr}})$ 
and consider 
{\scriptsize
\begin{align*}
P-(r)m_{4_{nr}} & = (r-n+1)m_{1_{nr}}+n\cdot m_{2_{nr}}+m_{3_{nr}}-m_{4_{nr}}-(r)m_{4_{nr}}\\
=& (r-n+1)s+n(s+3)+s+3n+1-(s+3n+2)-(r)(s+3n+2)\\
=&s(r-n+1+n-r)+3n-1-(r)(3n+2)\\
=&s+3n-1-r(3n-2)=3n+2>0.
\end{align*}
}
Therefore, we get $ P>rm_{4_{nr}}$ and $p+s=(r+1)m_{4_{nr}}$. Thus, $\gamma_{p}=s$.

Suppose $Q=(n-2)m_{1_{nr}}+n \cdot m_{2_{nr}}+2m_{3_{nr}}-m_{4_{nr}} \in \mathrm{PF}(\Gamma_{B_{nr}})$. We 
have,
{\scriptsize 
\begin{align*}
Q-(2n-2)m_{4_{nr}} & = (n-2)m_{1_{nr}}+n \cdot m_{2_{nr}}+2m_{3_{nr}}-m_{4_{nr}}-(2n-2)m_{4_{nr}}\\
&=(n-2)s+n (s+3)+2(s+3n+1)-(s+3n+2)-(2n-2)(s+3n+2)\\
&=s(n-2+n+2-1-2n+2)+3n+6n+2-3n-3-6n-4n+6n+4\\
&=s+2n+3>0.
\end{align*}
}
Therefore, $ Q>(2n-2)m_{4_{nr}}$ and $Q+(6n^{2}-5n-2)=(2n-1)m_{4_{nr}}$, hence $\gamma_{Q}=6n^{2}-5n-2$. 

Finally, we consider $R'=(r-2n+2)m_{1_{nr}}+2n\cdot m_{2_{nr}}-m_{4_{nr}}\in \mathrm{PF}(\Gamma_{B_{nr}})$. Then 
{\scriptsize
\begin{align*}
R'-rm_{4_{nr}} & = (r-2n+2)s+2n(s+3)-(s+3n+2)-r(s+3n+2)\\
= & s(r-2n+2+2n-1-r)+6n-3n-3-r(3n+2)\\
= & s+3n-3-r(3n+2)=3n>0.
\end{align*}
}
Therefore, $ R'>rm_{4_{nr}}$ and $R'+(s+1)=(r+1)m_{4_{nr}}$, thus $\gamma_{R'}=s+1$. 

It follows that $P,Q,R'$ satisfy property ($*$) in Theorem \ref{generator}, hence 
$D_{2}=\displaystyle\bigcup_{i=2}^{7}B_{i}$. \qed

\begin{remark}{\rm
Note that $| \mathrm{PF}(\mathrm{pr}_{1}(\overline{\Gamma}_{B_{nr}})|=| \mathrm{PF}(\Gamma_{B_{nr}})|=3n+2 $ and $| \mathrm{PF}(\mathrm{pr}_{2}(\overline{\Gamma}_{B_{nr}}))|=1$. Therefore, by Theorem \ref{Poincare-1}, $P_{\mathrm{Der}_{\mathbb{K}}(R)}^{R}(z)= 1+(3n+3)P_{\mathbb{K}}^{R}(z)$. 
\medskip

The problem of finding the Poincar\'{e} series of derivation modules 
of the projective closure $\mathrm{Der}_{\mathbb{K}}(\overline{C(n_{1},\ldots,n_{e})})$ over $\overline{C(n_{1},\ldots,n_{e})}$ is reduced to finding the Poincar\'{e} series of derivation modules of $\mathbb{K}$ over $\overline{C(n_{1},\ldots,n_{e})}$. This leads us to 
the following question:
\medskip

\noindent\textbf{Question:} What is the Poincar\'{e} series of the field $\mathbb{K}$ over the projective closure $\overline{C(n_{1},\ldots,n_{e})}$?
}
\end{remark}

\end{document}